\documentclass{amsart}%
\usepackage{amsthm}
\usepackage{amsmath}
\usepackage{amsfonts}
\usepackage{amssymb}
\usepackage{dsfont}
\usepackage{esint}
\usepackage{graphicx}
\usepackage{color}
\usepackage{url}
\usepackage[latin1]{inputenc}
\usepackage[T1]{fontenc}
\usepackage[english]{babel}%
\numberwithin{equation}{section}

\usepackage{amsaddr}

\providecommand{\U}[1]{\protect\rule{.1in}{.1in}}

\newtheorem{definition}{Definition}[section]

\newtheorem{question}{Question}[section]

\newcommand{\rst}[1]{\ensuremath{\raise-1.0ex\hbox{\large{$\vert_{#1}$}}}}

\begin{document}

\title[Existence for a class of viscous fluid problems]{A note about existence for a class of viscous fluid problems}

\date{\textbf{March 30, 2012}}
\author[H.B. de Oliveira]{Hermenegildo Borges de Oliveira$^{\ast,\ast\ast}$}

\email{holivei@ualg.pt}

\address{$^{\ast}$FCT - Universidade do Algarve and $^{\ast\ast}$CMAF - Universidade de Lisboa, Portugal.}

\begin{abstract}
In this work the existence of weak solutions for a class of non-Newtonian viscous fluid problems is analyzed.
The problem is modeled by the steady case of the generalized Navier-Stokes equations, where the exponent $q$ that characterizes the flow depends on the space variable: $q=q(\mathbf{x})$. For the associated boundary-value problem we show that, in some situations, the log-Hölder continuity condition on $q$ can be dropped and the result of the existence of weak solutions still remain valid for any variable exponent $q\geq\alpha>\frac{2N}{N+2}$, where $\alpha=\mathrm{ess}\inf q$.



\end{abstract}

\maketitle

\selectlanguage{english}





\section{Introduction}\label{Sect-Int}

In this article we study the steady motion of an incompressible and homogeneous viscous fluid in a bounded domain $\Omega\subset\mathds{R}^N$, $N\geq 2$, with the boundary denoted by $\partial\Omega$. We assume the motion is described by the following boundary-value problem for the generalized Navier-Stokes equations:
\begin{equation}\label{geq1-inc}
\mathrm{div}\,\mathbf{u}=0\quad\mbox{in}\quad \Omega;
\end{equation}
\begin{equation}\label{geq1-vel}
\mathbf{div}(\mathbf{u}\otimes\mathbf{u})
=\mathbf{f}-\mathbf{\nabla}p+\mathbf{div}\,\mathbf{S} \quad\mbox{in}\quad \Omega;
\end{equation}
\begin{equation}\label{geq1-bc}
\mathbf{u}=\mathbf{0} \qquad\mbox{on}\quad \partial\Omega.
\end{equation}
Here, $\mathbf{u}$ is the velocity field, $p$ stands for the pressure divided by
the constant density and $\mathbf{f}$ is the external forces field. %
We assume the extra stress tensor $\mathbf{S}$ has a variable $q$-structure in the following sense:
\begin{itemize}
\item[(A)]
$\mathbf{S}:\Omega\times\mathds{M}^N_{\mathrm{sym}}\to\mathds{M}^N_{\mathrm{sym}}$ is a Charathéodory function;
\item[(B)] $|\mathbf{S}(\mathbf{x},\mathbf{A})|\leq C|\mathbf{A}|^{q(\mathbf{x})-1}$ for all
$\mathbf{A}$ in $\mathds{M}^N_{\mathrm{sym}}$ and a.a. $\mathbf{x}$ in $\Omega$;
\item[(C)] $\mathbf{S}(\mathbf{x},\mathbf{A}):\mathbf{A}\geq C|\mathbf{A}|^{q(\mathbf{x})}$ for all
$\mathbf{A}$ in $\mathds{M}^N_{\mathrm{sym}}$ and a.a. $\mathbf{x}$ in $\Omega$;
\item[(D)] $\left(\mathbf{S}(\mathbf{x},\mathbf{A})-\mathbf{S}(\mathbf{x},\mathbf{B})\right):(\mathbf{A}-\mathbf{B})>0$ for all
$\mathbf{A}\not=\mathbf{B}$ in $\mathds{M}^N_{\mathrm{sym}}$ and a.a. $\mathbf{x}$ in $\Omega$.
\end{itemize}
Here, $\mathds{M}^N_{\mathrm{sym}}$ is the vector space of all symmetric $N\times N$ matrices, which is equipped with the scalar product $\mathbf{A}:\mathbf{B}$ and
norm $|\mathbf{A}|=\sqrt{\mathbf{A}:\mathbf{A}}$. %

The existence of weak solutions to the problem (\ref{geq1-inc})-(\ref{geq1-bc}) with a constant $q$-structure was established by \cite{L-1967} and \cite{Lions-1969} for $q\geq\frac{3N}{N+2}$, by \cite{FMS-1997} and \cite{R-1997} for $q>\frac{2N}{N+1}$ and, finally and again, by \cite{FMS-2003} for $q>\frac{2N}{N+2}$. %
These results were obtained in the class
\begin{equation}\label{V-q=const}
\mathbf{V}_q:=\mbox{closure of $\mathcal{V}$ in $\mathbf{W}^{1,q}(\Omega)$},\quad\mbox{where}\quad\mathcal{V}:=\{\mathbf{v}\in\mathbf{C}^{\infty}_0(\Omega):\mathrm{div\,}\mathbf{v}=0\}\,.
\end{equation}
The proofs in \cite{L-1967,Lions-1969} use the theory of monotone operators together with compactness arguments, whereas in \cite{FMS-1997,R-1997} and \cite{FMS-2003} are used, in addition, the $L^{\infty}$  and the Lipschitz-truncation methods, respectively. %
Each one of these results improves the previous one in the sense that the convective term $\mathbf{u}\otimes\mathbf{u}:\mathbf{D}(\mathbf{\varphi})$ is in $\mathbf{L}^1(\Omega)$ for an increasingly smaller lower limit of $q$.

The mathematical analysis of the problem (\ref{geq1-inc})-(\ref{geq1-bc}), with the deviatoric stress tensor satisfying to (A)-(D) with a variable $q$-structure, must be done in the context of Orlicz spaces. These spaces resemble many of the aspects of classical Lebesgue and Sobolev spaces, but there are some important differences which must be pointed out (see Section~\ref{Sect-Not}). %
Existence results for the problem (\ref{geq1-inc})-(\ref{geq1-bc}), with the deviatoric stress tensor satisfying to (A)-(D) with a variable $q$-structure, are due to \cite{R-2000}, \cite{H-2011} and \cite{DMS-2008} and were obtained in the class
\begin{equation}\label{Wq-nolog}
\mathbf{W}_q(\Omega):=\mbox{closure of $\mathcal{V}$ in the $\|\mathbf{D(v)}\|_{\mathbf{L}^{q(\cdot)}(\Omega)}$--\ norm}\,,
\end{equation}
where $q\in\mathcal{P}(\Omega)$, the set of all measurable functions $q:\Omega\to [1,\infty]$, satisfies to
\begin{equation}\label{ap1}
\displaystyle
1<\alpha:=\mathrm{ess}\inf_{\hspace{-0.5cm}\mathbf{x}\in \Omega}q(\mathbf{x})\leq q(\mathbf{x})\leq \mathrm{ess}\sup_{\hspace{-0.5cm}\mathbf{x}\in \Omega}q(\mathbf{x}):=\beta<\infty.
\end{equation}
The proofs in~\cite{R-2000} and in \cite{H-2011} are valid for $\alpha>\frac{3N}{N+2}$ and $\alpha>\frac{2N}{N+1}$, respectively. Moreover they follow the same approach of \cite{L-1967,Lions-1969} and \cite{FMS-1997,R-1997}, respectively, and use the fact that $\mathbf{W}_q(\Omega)$ is continuously imbedded into $\mathbf{V}_{\alpha}$. The proof of~\cite{DMS-2008} is valid for $\alpha>\frac{2N}{N+2}$, provided the variable exponent $q$ is globally log-Hölder continuous in the sense of (\ref{glob-log-H}) below. %
The proof here follows the same approach of the result for constant $q$ in \cite{DMS-2008} and uses results on Lipschitz truncations of functions in Orlicz-Sobolev spaces performed still in~\cite{DMS-2008}. See also \cite{ADO-PNDEA,AR-2006,R-1997} for concrete fluid models with a variable $q$-structure.

Our goal in this work is to show that the log-Hölder continuity condition (\ref{glob-log-H}) is not necessary to show the existence of weak solutions to the problem (\ref{geq1-inc})-(\ref{geq1-bc}) with the deviatoric stress tensor satisfying to (A)-(D) with a variable $q$-structure. %
 As one can sees in the proof of \cite[Theorem 5.1]{DMS-2008}, assumption (\ref{glob-log-H}) is fundamental to achieve the existence result by the method proposed there. %
Firstly, we shall seek for a different condition that assures the existence of weak solutions for this problem in the case of $\alpha>\frac{2N}{N+2}$. At the end, we shall give an example to which neither this new condition nor (\ref{glob-log-H}) are needed.

\section{Weak Formulation}\label{Sect-Not}
\numberwithin{equation}{section}

The notation used in this work is largely standard in Mathematical Fluid Mechanics (see \emph{e.g.} \cite{Lions-1969}).
In this article, the notations $\Omega$ or $\omega$ stand always for a domain, \emph{i.e.}, a connected open subset of $\mathds{R}^N$, $N\geq 1$.
Given $k\in\mathds{N}$, we denote by $\mathrm{C}^{k}(\Omega)$ the space of all $k$-differentiable functions in $\Omega$. %
By $\mathrm{C}^{\infty}_0(\Omega)$ we denote the space of all infinity-differentiable functions with compact support in $\Omega$. %
In the context of distributions, the space $\mathrm{C}^{\infty}_0(\Omega)$ is denoted by $\mathcal{D}(\Omega)$ instead. %
The space of distributions over $\mathcal{D}(\Omega)$ is denoted by $\mathcal{D}'(\Omega)$. %
If $\mathrm{X}$ is a generic Banach space, its dual space is denoted by $\mathrm{X}'$.
Let $1\leq q\leq \infty$ and $\Omega\subset\mathds{R}^N$, with $N\geq 1$, be a domain. %
We use the classical Lebesgue spaces $\mathrm{L}^q(\Omega)$, whose norm is
denoted by $\|\cdot\|_{\mathrm{L}^q(\Omega)}$. %
For any nonnegative $k$,
$\mathrm{W}^{k,q}(\Omega)$ denotes the Sobolev space of all
functions $u\in\mathrm{L}^q(\Omega)$ such
that the weak derivatives $\mathrm{D}^{\alpha}u$ exist, in the generalized sense, and are in
$\mathrm{L}^q(\Omega)$ for any multi-index $\alpha$ such that
$0\leq |\alpha|\leq k$.
In particular, $\mathrm{W}^{1,\infty}(\Omega)$ stands for the space of Lipschitz functions. %
The norm in $\mathrm{W}^{k,q}(\Omega)$ is denoted by
$\|\cdot\|_{\mathrm{W}^{k,q}(\Omega)}$.
We define $\mathrm{W}^{k,q}_0(\Omega)$ as the closure of $\mathrm{C}^{\infty}_0(\Omega)$ in $\mathrm{W}^{k,q}(\Omega)$.
For the dual space of $\mathrm{W}^{k,q}_0(\Omega)$, we use the identity $(\mathrm{W}^{k,q}_0(\Omega))'=\mathrm{W}^{-k,q'}(\Omega)$, up to an isometric isomorphism. %
Vectors and vector spaces will be denoted by boldface letters.

We denote by $\mathcal{P}(\Omega)$ the set of all measurable functions $q:\Omega\to [1,\infty]$ and define
$$\displaystyle q^{-}:=\mathrm{ess}\inf_{\hspace{-0.5cm}x\in \Omega}q(x),\quad q^{+}:=\mathrm{ess}\sup_{\hspace{-0.5cm}x\in \Omega}q(x).$$
Given $q\in\mathcal{P}(\Omega)$, we denote by $\mathrm{L}^{q(\cdot)}(\Omega)$ the space of all measurable functions $f$ in $\Omega$ such that its semimodular is finite:
\begin{equation}\label{ap3}
A_{q(\cdot)}(f):=\int_{\Omega}|f(x)\mathbf{|}^{q(x)}d\,x<\infty.
\end{equation}
The space $\mathrm{L}^{q(\cdot)}(\Omega)$ is called Orlicz space and is also known by Lebesgue space with variable exponent. %
Equipped with the norm
\begin{equation}\label{ap4}
\left\Vert f\right\Vert_{\mathrm{L}^{q(\cdot)}(\Omega)}:=
\inf\left\{\lambda>0 : A_{q(\cdot)}\left(\frac{f}{\lambda}\right) \leq
1\right\},
\end{equation}
$\mathrm{L}^{q(\cdot)}(\Omega)$ becomes a Banach space. %
If $q^{+}<\infty$, $\mathrm{L}^{q(\cdot)}(\Omega)$ is separable and the space $\mathrm{C}^{\infty}_0(\Omega)$ is dense in $\mathrm{L}^{q(\cdot)}(\Omega)$.
Moreover, if $1<q^{-}\leq q^{+}<\infty$, $\mathrm{L}^{q(\cdot)}(\Omega)$ is reflexive. %
One problem in Orlicz spaces is the relation between the semimodular (\ref{ap3}) and the norm (\ref{ap4}). %
If (\ref{ap1}) is satisfied, one can shows that
\begin{equation*}\label{ap4-56}
\|f\|_{\mathrm{L}^{q(\cdot)}(\Omega)}^{q^{-}}-1\leq A_{q(\cdot)}(f)\leq
\|f\|_{\mathrm{L}^{q(\cdot)}(\Omega)}^{q^{+}}+1\,.
\end{equation*}
In Orlicz spaces, there holds a version of H\"{o}lder's inequality, called generalized H\"{o}lder's inequality. %
Given $q\in\mathcal{P}(\Omega)$, the Orlicz-Sobolev space $W^{1,q(\cdot)}(\Omega)$ is defined as the set of all functions $f\in\mathrm{L}^{q(\cdot)}(\Omega)$ such that $\mathrm{D}^{\alpha}f\in\mathrm{\mathrm{L}}^{q(\cdot)}(\Omega)$ for any multi-index $\alpha$ such that $0\leq |\alpha|\leq 1$.
In $W^{1,q(\cdot)}(\Omega)$ is defined a semimodular and the correspondent induced norm analogously as in (\ref{ap3})-(\ref{ap4}). %
For this norm, $W^{1,q(\cdot)}(\Omega)$ is a Banach space, which becomes separable and reflexive in the same conditions as $\mathrm{L}^{q(\cdot)}(\Omega)$. %
The Orlicz-Sobolev space with zero boundary values is defined by:
$$W^{1,q(\cdot)}_0(\Omega):=\overline{\left\{f\in W^{1,q(\cdot)}(\Omega): \mathrm{supp}\ f\subset\subset \Omega\right\}}^{\ \|\cdot\|_{W^{1,q(\cdot)}(\Omega)}}\,.$$
In contrast to the case of classical Sobolev spaces, the set $\mathrm{C}_0^{\infty}(\Omega)$ is not necessarily dense in $\mathrm{W}^{1,q(\cdot)}_0(\Omega)$ -- the closure of $\mathrm{C}_0^{\infty}(\Omega)$ in $\mathrm{W}^{1,q(\cdot)}(\Omega)$ is strictly contained in $\mathrm{W}^{1,q(\cdot)}_0(\Omega)$. %
The equality holds only if $q$ is globally log-Hölder continuous, \emph{i.e.}, if exist positive constants
$C_1$, $C_2$ and $q_{\infty}$ such that
\begin{equation}\label{glob-log-H}
\left|q(\mathbf{x})-q(\mathbf{y})\right|\leq \frac{C_1}{\ln(e+1/|\mathbf{x}-\mathbf{y}|)},\quad
\left|q(\mathbf{x})-q_{\infty}\right|\leq \frac{C_2}{\ln(e+|\mathbf{x}|)}\quad\forall\ \mathbf{x},\ \mathbf{y}\in \Omega.
\end{equation}
See the monograph \cite{DHHR} for a thorough analysis on Orlicz and Orlicz-Sobolev spaces.

In order to introduce the notion of weak solutions we shall consider in this work, let us recall the well-known function spaces of Mathematical Fluid Mechanics defined at (\ref{V-q=const}). Due to the presence of the variable exponent $q(\cdot)$ in the structure of the tensor $\mathbf{S}$, we need to consider the weak solutions to the problem (\ref{geq1-inc})-(\ref{geq1-bc}) in some Orlicz-Sobolev space. %
Since the set $\mathbf{C}_0^{\infty}(\Omega)$ is not necessarily dense in $\mathbf{W}^{1,q(\cdot)}_0(\Omega)$, we define the analogue of $\mathbf{V}_q$ by (\ref{Wq-nolog}). %
It is a easy task to verify the space $\mathbf{W}_q(\Omega)$ satisfies to the following continuous imbeddings:
\begin{equation*}
\mathbf{V}_{\beta}\hookrightarrow\mathbf{W}_q(\Omega)\hookrightarrow\mathbf{V}_{\alpha}\,.
\end{equation*}
Moreover, $\mathbf{W}_q(\Omega)$ is a closed subspace of $\mathbf{V}_{\alpha}$ and therefore it is
a reflexive and separable Banach space for the norm
\begin{equation*}
\|\mathbf{v}\|_{\mathbf{W}_q(\Omega)}:=
\|\mathbf{D(v)}\|_{\mathbf{L}^{q(\cdot)}(\Omega)}.
\end{equation*}

\begin{definition}\label{weak-sol-vq}
Let $\Omega$ be a bounded domain of $\mathds{R}^N$, with $N \geq 2$, and let $q\in\mathcal{P}(\Omega)$ be a variable exponent satisfying to (\ref{ap1}). %
Let also $\mathbf{f}\in\mathbf{L}^{1}(\Omega)$ and assume that conditions (A)-(D) are fulfilled with a variable exponent $q$.
A vector field $\mathbf{u}$ is a weak solution to the problem (\ref{geq1-inc})-(\ref{geq1-bc}), if:
\begin{enumerate}
\item $\mathbf{\mathbf{u}}\in  \mathbf{W}_q(\Omega)$;
\item For every $\varphi\in\mathbf{C}^{\infty}_0(\Omega)$, with $\mathrm{div\,}\varphi=0$,
\begin{equation*}
\int_{\Omega}\left(\mathbf{S}(\mathbf{D}(\mathbf{u}))-\mathbf{u}\otimes\mathbf{u}\right):\mathbf{D}(\mathbf{\varphi})\,d\mathbf{x}
=\int_{\Omega}\mathbf{f}\cdot\varphi\,d\,\mathbf{x}.
\end{equation*}
\end{enumerate}
\end{definition}

\noindent The main goal of this work is to seek for the condition(s) we have to impose in the problem (\ref{geq1-inc})-(\ref{geq1-bc}) that assure(s) the existence of weak solutions to this problem in the sense of Definition~\ref{weak-sol-vq} and without any further restriction on $q$ besides
(\ref{ap1}) above and  (\ref{h-f})-(\ref{al-TE}) below.

\begin{question}\label{th-exst-ws-pp}
Let $\Omega$ be a bounded domain in $\mathds{R}^{N}$, $N\geq 2$. %
Assume that conditions (A)-(D) are fulfilled with a variable exponent $q\in\mathcal{P}(\Omega)$ satisfying to (\ref{ap1}), and
\begin{equation}\label{h-f}
\mathbf{f}=-\mathbf{div}\,\mathbf{F},\quad \mathbf{F}\in\mathds{M}^N_{\mathrm{sym}}\,,
\quad \mathbf{F}\in\mathbf{L}^{q'(\cdot)}(\Omega),
\end{equation}
\begin{equation}\label{al-TE}
\frac{2N}{N+2}<\alpha\leq\beta<\infty\,.
\end{equation}
Is it possible to find a distinct condition from the log-Hölder continuity property (\ref{glob-log-H}) that assures the existence of a weak solution to the problem (\ref{geq1-inc})-(\ref{geq1-bc}) in the sense of Definition~\ref{weak-sol-vq}?
\end{question}

\noindent The answer to Question~\ref{th-exst-ws-pp} will be the aim of next sections. For that, we shall prove an existence result for the problem (\ref{geq1-inc})-(\ref{geq1-bc}) under the conditions stated in Question~\ref{th-exst-ws-pp}. We will see that the validity of such an existence result will demand a new and different condition.

\section{The regularized problem}\label{Sect-Exist-RP}

Let $\Phi\in\mathrm{C}^{\infty}([0,\infty))$ be a non-increasing function such that $0\leq\Phi\leq 1$ in $[0,\infty)$, $\Phi\equiv 1$ in $[0,1]$, $\Phi\equiv 0$ in $[2,\infty)$ and $0\leq -\Phi'\leq 2$. %
For $\epsilon>0$, we set
\begin{equation}\label{Phi-e}
\Phi_{\epsilon}(s):=\Phi(\epsilon s),\quad s\in[0,\infty).
\end{equation}
We consider the following regularized problem:
\begin{equation}\label{eq2-inc-e}
\mathrm{div}\,\mathbf{u}_{\epsilon}=0\quad\mbox{in}\quad \Omega,
\end{equation}
\begin{equation}\label{eq2-vel-e}
\mathbf{div}(\mathbf{u}_{\epsilon}\otimes\mathbf{u}_{\epsilon}\Phi_{\epsilon}(|\mathbf{u}_{\epsilon}|))
=\mathbf{f}-\mathbf{\nabla}p_{\epsilon}+
\mathbf{div}\left(\mathbf{S}(\mathbf{D}(\mathbf{u}_{\epsilon}))+\epsilon|\mathbf{D}(\mathbf{u}_{\epsilon})|^{\beta-2}\mathbf{D}(\mathbf{u}_{\epsilon})\right)
\quad\mbox{in}\quad \Omega,
\end{equation}
\begin{equation}\label{eq1-bc-u-e}
\mathbf{u}_{\epsilon}=\mathbf{0}\qquad\mbox{on}\quad \partial\Omega.
\end{equation}
A vector function $\mathbf{u}_{\epsilon}\in\mathbf{V}_{\beta}$ is a weak solution to the problem (\ref{eq2-inc-e})-(\ref{eq1-bc-u-e}), if
\begin{equation}\label{eq-ws-reg}
\int_{\Omega}\left[\mathbf{S}(\mathbf{D}(\mathbf{u}_{\epsilon}))+\epsilon|\mathbf{D}(\mathbf{u}_{\epsilon})|^{\beta-2}\mathbf{D}(\mathbf{u}_{\epsilon})-\mathbf{u}_{\epsilon}\otimes\mathbf{u}_{\epsilon}\Phi_{\epsilon}(|\mathbf{u}_{\epsilon}|)\right]:\mathbf{D}(\mathbf{\varphi})\,d\mathbf{x}
=\int_{\Omega}\mathbf{F}:\mathbf{D}(\varphi)\,d\,\mathbf{x}
\end{equation}
for all $\varphi\in\mathcal{V}$. %
Under the assumptions stated in Question~\ref{th-exst-ws-pp}, it can be proved that, for each $\epsilon>0$, there exists a weak solution $\mathbf{u}_{\epsilon}\in\mathbf{V}_{\beta}$ to the problem (\ref{eq2-inc-e})-(\ref{eq1-bc-u-e}). %
The proof of this result is based on the Schauder fixed point  theorem. The map construction is done by putting the convective term on the  right hand side and by solving a nonlinear equation via the monotone operator theory. %
Moreover, it can be proved that every weak solution satisfies to the following energy equality:
\begin{equation}\label{e-equality-qr}
\int_{\Omega}\mathbf{S}(\mathbf{D}(\mathbf{u}_{\epsilon})):\mathbf{D}(\mathbf{u}_{\epsilon})d\mathbf{x}+
\epsilon\int_{\Omega}|\mathbf{D}(\mathbf{u}_{\epsilon})|^{\beta}d\mathbf{x}
=\int_{\Omega}\mathbf{F}:\mathbf{D}(\mathbf{u}_{\epsilon})d\mathbf{x}.
\end{equation}

\noindent Now, let $\mathbf{u}_{\epsilon}\in\mathbf{V}_{\beta}$ be a weak solution to the problem (\ref{eq2-inc-e})-(\ref{eq1-bc-u-e}). %
From  (\ref{e-equality-qr}) we can prove that
\begin{equation}\label{e-inequality-qr}
\int_{\Omega}|\mathbf{D}(\mathbf{u}_{\epsilon})|^{q(\mathbf{x})}d\mathbf{x}+
\epsilon\int_{\Omega}|\mathbf{D}(\mathbf{u}_{\epsilon})|^{\beta}d\mathbf{x}\leq C,
\end{equation}
where, by the assumption (\ref{h-f}), $C$ is a positive constant and, very important, does not depend on $\epsilon$. %
Then we can prove from (\ref{e-inequality-qr}) that
\begin{equation}\label{est-q}
\|\mathbf{D}(\mathbf{u}_{\epsilon})\|_{\mathbf{L}^{q(\cdot)}(\Omega)}\leq C,
\end{equation}
\begin{equation}\label{est-inf-gam}
\|\mathbf{u}_{\epsilon}\|_{\mathbf{V}_{\alpha}}\leq C,
\end{equation}
\begin{equation}\label{est-q'}
\|\mathbf{S}(\mathbf{D}(\mathbf{u}_{\epsilon}))\|_{\mathbf{L}^{q'(\cdot)}(\Omega)}\leq C,
\end{equation}
\begin{equation}\label{est-gam'}
\|\mathbf{S}(\mathbf{D}(\mathbf{u}_{\epsilon}))\|_{\mathbf{L}^{\beta'}(\Omega)}\leq C.
\end{equation}
On the other hand, by using (\ref{est-inf-gam}) and Sobolev's inequality, we have
\begin{equation}\label{est-gam(N+2)/N}
\|\mathbf{u}_{\epsilon}\|_{\mathbf{L}^{\alpha^{\ast}}(\Omega)}\leq C,
\end{equation}
where $\alpha^{\ast}$ denotes the Sobolev conjugate of $\alpha$. %
As a consequence of (\ref{est-gam(N+2)/N}) and (\ref{Phi-e}),
\begin{equation}\label{est-Phiuxu}
\|\mathbf{u}_{\epsilon}\otimes\mathbf{u}_{\epsilon}\Phi_{\epsilon}(|\mathbf{u}_{\epsilon}|)\|_{\mathbf{L}^{\frac{\alpha^{\ast}}{2}}(\Omega)}\leq C.
\end{equation}
From (\ref{est-inf-gam}), (\ref{est-gam'}) and (\ref{est-Phiuxu}), there exists a sequence of positive numbers $\epsilon_m$ such that
$\epsilon_m\to 0$, as $m\to \infty$, and
\begin{equation}\label{convg-La}
\mbox{$\mathbf{u}_{\epsilon_m}\to \mathbf{u}$\quad weakly in $\mathbf{V}_{\alpha}$,\quad as $m\to\infty$,}
\end{equation}
\begin{equation}\label{convg-S-b'}
\mbox{$\mathbf{S}(\mathbf{D}(\mathbf{u}_{\epsilon_m}))\to \mathbf{S}$\quad weakly in $\mathbf{L}^{\beta'}(\Omega)$,\quad as $m\to\infty$,}
\end{equation}
\begin{equation}\label{convg-Phiuxu}
\mbox{$\mathbf{u}_{\epsilon_m}\otimes\mathbf{u}_{\epsilon_m}\Phi_{\epsilon_m}(|\mathbf{u}_{\epsilon_m}|)\to \mathbf{G}$\quad weakly in $\mathbf{L}^{\frac{\alpha^{\ast}}{2}}(\Omega)$,\quad as $m\to\infty$},
\end{equation}
\begin{equation}\label{conv-Re-0}
\epsilon_{m}|\mathbf{D}(\mathbf{u}_{\epsilon_m})|^{q-2}\mathbf{D}(\mathbf{u}_{\epsilon_m})\to 0\quad\mbox{weakly in}\quad\mathbf{L}^{\beta'}(\Omega),\quad\mbox{as}\quad m\to\infty.
\end{equation}

Now we observe that, due to (\ref{convg-La}), the application of Sobolev's compact imbedding theorem implies
\begin{equation}\label{convg-s-g}
\mbox{$\mathbf{u}_{\epsilon_m}\to \mathbf{u}$\quad strongly in $\mathbf{L}^{\gamma}(\Omega)$,\quad as $m\to\infty$,
\quad
for any $\gamma:1\leq\gamma<\alpha^{\ast}$.
}
\end{equation}
Since $2<\alpha^{\ast}$, it follows from (\ref{convg-s-g}) that
\begin{equation}\label{convg-s-2}
\mbox{$\mathbf{u}_{\epsilon_m}\to \mathbf{u}$\quad strongly in $\mathbf{L}^{2}(\Omega)$,\quad as $m\to\infty$.
}
\end{equation}
Using (\ref{Phi-e}) and (\ref{convg-s-2}), we can prove that
\begin{equation}\label{limit-eq-G}
\mbox{$\mathbf{u}_{\epsilon_m}\otimes\mathbf{u}_{\epsilon_m}\Phi_{\epsilon_m}(|\mathbf{u}_{\epsilon_m}|)\to \mathbf{u}\otimes\mathbf{u}$\quad strongly in $\mathbf{L}^{1}(\Omega)$,\quad as $m\to\infty$.
}
\end{equation}
Then gathering the information of (\ref{convg-Phiuxu}) and (\ref{limit-eq-G}), we see that $\mathbf{G}=\mathbf{u}\otimes\mathbf{u}$.

Finally, using the convergence results (\ref{convg-La})-(\ref{conv-Re-0}) and observing (\ref{limit-eq-G}), we can pass to the limit $m\to\infty$ in the following integral identity, which results from (\ref{eq-ws-reg}),
\begin{equation}\label{eq-wsm-reg}
\int_{\Omega}\left[\mathbf{S}(\mathbf{D}(\mathbf{u}_{\epsilon_m}))+\epsilon_m|\mathbf{D}(\mathbf{u}_{\epsilon_m})|^{\beta-2}\mathbf{D}(\mathbf{u}_{\epsilon_m})-\mathbf{u}_{\epsilon_m}\otimes\mathbf{u}_{\epsilon_m}\Phi_{\epsilon_m}(|\mathbf{u}_{\epsilon_m}|)-\mathbf{F}\right]:\mathbf{D}(\mathbf{\varphi})\,d\mathbf{x}
=0,
\end{equation}
valid for all $\mathbf{\varphi}\in\mathcal{V}$, to obtain
\begin{equation}\label{limit-eq-SH}
\int_{\Omega}(\mathbf{S}-\mathbf{u}\otimes\mathbf{u}-\mathbf{F}):\mathbf{D}(\mathbf{\varphi})\,d\mathbf{x}=0\quad
\forall\ \varphi\in\mathcal{V}. %
\end{equation}

\section{Decomposition of the pressure.}\label{Sect-Dec-p}

Since we shall use test functions which are not divergence free, we first have to determine the approximative pressure from the weak formulation
(\ref{eq-wsm-reg}). %
First, let $\omega'$ be a fixed but arbitrary open bounded subset of $\Omega$ such that
\begin{equation}\label{om-l}
\omega'\subset\subset\Omega\quad\mbox{and}\quad\partial\omega'\ \mbox{is Lipschitz}, %
\end{equation}
where $\omega'\subset\subset\Omega$ means that $\omega'$ is compactly contained in $\Omega$,
and let us set
\begin{equation}\label{Q-em}
\mathbf{Q}_{\epsilon_m}:=\mathbf{S}(\mathbf{D}(\mathbf{u}_{\epsilon_m}))+\epsilon_m|\mathbf{D}(\mathbf{u}_{\epsilon_m})|^{\beta-2}\mathbf{D}(\mathbf{u}_{\epsilon_m})-\mathbf{u}_{\epsilon_m}\otimes\mathbf{u}_{\epsilon_m}\Phi_{\epsilon_m}(|\mathbf{u}_{\epsilon_m}|)-\mathbf{F}.
\end{equation}
Using assumption (\ref{h-f}) and the results (\ref{est-gam'}), (\ref{est-Phiuxu}) and (\ref{conv-Re-0}), we can prove that
\begin{equation}\label{Q-in-Lr}
\mathbf{Q}_{\epsilon_m}\in\mathbf{L}^r(\Omega),\quad\mbox{where}\quad
 1<r\leq r_0:=\min\left\{\beta',\frac{\alpha^{\ast}}{2}\right\}.
\end{equation}
Note that $r_0=\min\{\beta',\frac{N\alpha}{2(N-\alpha)}\}$ if $\alpha<N$ and $r_0=\beta'$ if $N\geq\alpha$. %
Then we define a linear functional
\begin{equation}\label{b1-lin-f-r}
\Pi_{\epsilon_m}:\mathbf{W}^{1,r'}_0(\omega')\rightarrow\mathbf{W}^{-1,r}(\omega')\,,
\end{equation}
\begin{equation}\label{b2-lin-f-r}
\langle\Pi_{\epsilon_m},\mathbf{\varphi}\rangle_{\mathbf{W}^{-1,r}(\omega')\times\mathbf{W}^{1,r'}_0(\omega')}:=
\int_{\omega'}\mathbf{Q}_{\epsilon_m}:\mathbf{D}(\mathbf{\varphi})\,d\mathbf{x}.
\end{equation}
Using (\ref{b1-lin-f-r})-(\ref{b2-lin-f-r}), we can prove, owing to (\ref{Q-in-Lr}),  that
\begin{equation}\label{b-lin-op-p}
\|\Pi_{\epsilon_m}\|_{(\mathbf{V}_{r'})'}\leq C\,,
\end{equation}
where $C$ is a positive constant independent of $m$. %
Note that here $\mathbf{V}_{r'}$ is taken over $\omega'$. %
Moreover, since $\mathcal{V}$ is dense in $\mathbf{V}_{r'}$, we can see that (\ref{eq-wsm-reg}), (\ref{Q-em}) and (\ref{b2-lin-f-r}) imply
\begin{equation}\label{dual=0}
\langle\Pi_{\epsilon_m},\mathbf{\varphi}\rangle_{(\mathbf{V}_{r'})'\times\mathbf{V}_{r'}}=0\quad\forall\
\mathbf{\varphi}\in\mathbf{V}_{r'}.
\end{equation}
By virtue of (\ref{b1-lin-f-r})-(\ref{dual=0}) and due to assumption (\ref{om-l}), we can apply a version of de Rham's Theorem to prove the existence of a unique function
\begin{equation}\label{app-pr}
p_{\epsilon_m}\in\mathbf{L}^{r'}(\omega'),\qquad\mbox{with}\quad\int_{\omega'}p_{\epsilon_m}d\mathbf{x}=0,
\end{equation}
such that
\begin{equation}\label{Rham-thm}
\langle\Pi_{\epsilon_m},\mathbf{\varphi}\rangle_{\mathbf{W}^{-1,r}(\omega')\times\mathbf{W}^{1,r'}_0(\omega')}=
\int_{\omega'}p_{\epsilon_m}\,\mathrm{div}\mathbf{\varphi}\,d\mathbf{x}\quad\forall\ \mathbf{\varphi}\in\mathbf{W}^{1,r'}_0(\omega')\,,
\end{equation}
\begin{equation}\label{b-pem}
\|p_{\epsilon_m}\|_{\mathbf{L}^{r'}(\omega')}\leq\|\Pi_{\epsilon_m}\|_{(\mathbf{V}_{r'})'}.
\end{equation}
Then, gathering the information of (\ref{eq-wsm-reg}), (\ref{Q-em}), (\ref{b2-lin-f-r}) and (\ref{Rham-thm}), we obtain
\begin{equation}\label{eq-wsm-reg-pem}
\begin{split}
 & \int_{\omega'}\mathbf{S}(\mathbf{D}(\mathbf{u}_{\epsilon_m})):\mathbf{D}(\mathbf{\varphi})\,d\mathbf{x}+
\epsilon_m\int_{\omega'}|\mathbf{D}(\mathbf{u}_{\epsilon_m})|^{\beta-2}\mathbf{D}(\mathbf{u}_{\epsilon_m}):\mathbf{D}(\mathbf{\varphi})\,d\mathbf{x}= \\
 & \int_{\omega'}\mathbf{F}:\mathbf{D}(\mathbf{\varphi})\,d\mathbf{x}+
\int_{\omega'}\mathbf{u}_{\epsilon_m}\otimes\mathbf{u}_{\epsilon_m}\Phi_{\epsilon_m}(|\mathbf{u}_{\epsilon_m}|):\mathbf{D}(\mathbf{\varphi})\,d\mathbf{x}+
\int_{\omega'}p_{\epsilon_m}\,\mathrm{div}\mathbf{\varphi}\,d\mathbf{x}
\end{split}
\end{equation}
for all $\mathbf{\varphi}\in\mathbf{W}^{1,r'}_0(\omega')$. %
On the other hand, due to (\ref{b-lin-op-p}) and (\ref{b-pem}) and by means of reflexivity, we get, passing to a subsequence, that
\begin{equation}\label{conv-p0-em}
p_{\epsilon_m}\to p_0\quad\mbox{weakly in}\quad\mathrm{L}^{r'}(\omega'),\quad\mbox{as}\quad m\to\infty.
\end{equation}
Next, passing to the limit $m\to\infty$ in the integral identity (\ref{eq-wsm-reg-pem}) by using the convergence results
(\ref{convg-S-b'}), (\ref{convg-Phiuxu}) together with (\ref{limit-eq-G}), using also (\ref{conv-Re-0}) and
(\ref{conv-p0-em}), we obtain
\begin{equation}\label{eq-limit-pem}
\begin{split}
 & \int_{\omega'}\left(\mathbf{S}-\mathbf{u}\otimes\mathbf{u}-\mathbf{F}\right):\mathbf{D}(\mathbf{\varphi})\,d\mathbf{x}=
\int_{\omega'}p_0\,\mathrm{div}\mathbf{\varphi}\,d\mathbf{x}
\end{split}
\end{equation}
for all $\mathbf{\varphi}\in\mathbf{W}^{1,r'}_0(\omega')$. %

Next, we shall decompose the pressure found in the first part of this section. %
With this in mind, let  $\omega$ be a fixed but arbitrary domain such that
\begin{equation}\label{om-C2}
\omega\subset\subset\omega'\subset\subset\Omega\quad\mbox{and}\quad\partial\omega\ \mbox{is $C^2$}. %
\end{equation}
To simplify the notation in the sequel, let us set
\begin{equation*}\label{set-A}
\mathrm{A}^s(\omega):=\{a\in\mathrm{L}^{s}(\omega):a=\triangle u,\quad u\in\mathrm{W}^{2,s}_0(\omega)\},\qquad  1<s<\infty. \\
\end{equation*}
Here we shall use some results due to \cite{Wolf-2007} that allow us to locally decompose the pressure. %
Applying \cite[Lemma 2.4]{Wolf-2007}, with $s=\beta'$ first and then with $s=\frac{\alpha^{\ast}}{2}$, and using (\ref{convg-S-b'}) and (\ref{conv-Re-0}) by one hand and (\ref{convg-Phiuxu}) and (\ref{limit-eq-G}) on the other, we can infer that exist unique functions
\begin{equation}\label{E1-p1}
p^1_{\epsilon_m}\in\mathrm{A}^{\beta'}(\omega),
\end{equation}
\begin{equation}\label{E1-p2}
p^2_{\epsilon_m}\in\mathrm{A}^{\frac{\alpha^{\ast}}{2}}(\omega)
\end{equation}
such that
\begin{equation}\label{dec-p1}
\int_{\omega}p^1_{\epsilon_m}\triangle\phi\,d\mathbf{x}=
\int_{\omega}\left(\mathbf{S}(\mathbf{D}(\mathbf{u}_{\epsilon_m}))+\epsilon_m|\mathbf{D}(\mathbf{u}_{\epsilon_m})|^{\beta-2}\mathbf{D}(\mathbf{u}_{\epsilon_m})- \mathbf{S}\right):\nabla^2\phi\,d\mathbf{x},
\end{equation}
\begin{equation}\label{dec-p2}
\int_{\omega}p^2_{\epsilon_m}\triangle\phi\,d\mathbf{x}=-
\int_{\omega}(\mathbf{u}_{\epsilon_m}\otimes\mathbf{u}_{\epsilon_m}\Phi_{\epsilon_m}(|\mathbf{u}_{\epsilon_m}|)-
\mathbf{u}\otimes\mathbf{u}):\nabla^2\phi\,d\mathbf{x}
\end{equation}
for all $\phi\in\mathrm{C}_0^{\infty}(\omega)$. %
Attending to (\ref{convg-S-b'}), (\ref{conv-Re-0}) and (\ref{dec-p1}) by one hand, and (\ref{convg-Phiuxu}), (\ref{limit-eq-G}) and (\ref{dec-p2}) on the other, a direct application of \cite[Lemma 2.3]{Wolf-2007}, with $s=\beta'$ and then with $s=\frac{\alpha^{\ast}}{2}$, yields
\begin{equation}\label{bd2-p1-em}
\|p^1_{\epsilon_m}\|_{\mathrm{L}^{\beta'}(\omega)}\leq C_1\|\mathbf{S}(\mathbf{D}(\mathbf{u}_{\epsilon_m}))-\mathbf{S}+\epsilon_m|\mathbf{D}(\mathbf{u}_{\epsilon_m})|^{\beta-2}\mathbf{D}(\mathbf{u}_{\epsilon_m})\|_{\mathbf{L}^{\beta'}(\omega)},
\end{equation}
\begin{equation}\label{bd2-p2-em}
\|p^2_{\epsilon_m}\|_{\mathrm{L}^{\frac{\alpha^{\ast}}{2}}(\omega)}\leq C_2\|\mathbf{u}_{\epsilon_m}\otimes\mathbf{u}_{\epsilon_m}\Phi_{\epsilon_m}(|\mathbf{u}_{\epsilon_m}|)-
\mathbf{u}\otimes\mathbf{u}\|_{\mathbf{L}^{\frac{\alpha^{\ast}}{2}}(\omega)}.
\end{equation}
where $C_1$ and $C_2$ are positive constants depending on $\beta'$, $\alpha^{\ast}$ and on $\omega$. %

Now, combining (\ref{eq-wsm-reg-pem}) and (\ref{eq-limit-pem}), and using the definition of the distributive derivative, we obtain
\begin{equation}\label{dist-der}
\begin{split}
   & \mathbf{div}\left(\mathbf{S}(\mathbf{D}(\mathbf{u}_{\epsilon_m}))-\mathbf{S}+\epsilon_m|\mathbf{D}(\mathbf{u}_{\epsilon_m})|^{\beta-2}\mathbf{D}(\mathbf{u}_{\epsilon_m})\right)-\\
   &\mathbf{div}\left(\mathbf{u}_{\epsilon_m}\otimes\mathbf{u}_{\epsilon_m}\Phi_{\epsilon_m}(|\mathbf{u}_{\epsilon_m}|)-\mathbf{u}\otimes\mathbf{u}\right)
   =\mathbf{\nabla}(p_{\epsilon_m}-p_0)
\end{split}
\qquad\mbox{in $\mathcal{D}'(\omega)$}.
\end{equation}
Then, testing (\ref{dist-der}) by $\nabla\phi$, with $\phi\in\mathrm{C}_0^{\infty}(\omega)$, integrating over $\omega$ and comparing the resulting equation with the one resulting from adding (\ref{dec-p1}) and (\ref{dec-p2}), we obtain
$$p_{\epsilon_m}-p_0=p^1_{\epsilon_m}+p^2_{\epsilon_m}\,.$$
Inserting this into (\ref{dist-der}), it follows that
\begin{equation}\label{dist2-der}
\begin{split}
   & \mathbf{div}\left(\mathbf{S}(\mathbf{D}(\mathbf{u}_{\epsilon_m}))-\mathbf{S}+\epsilon_m|\mathbf{D}(\mathbf{u}_{\epsilon_m})|^{\beta-2}\mathbf{D}(\mathbf{u}_{\epsilon_m})\right)-\\
   &\mathbf{div}\left(\mathbf{u}_{\epsilon_m}\otimes\mathbf{u}_{\epsilon_m}\Phi_{\epsilon_m}(|\mathbf{u}_{\epsilon_m}|)-\mathbf{u}\otimes\mathbf{u}\right)
   =\mathbf{\nabla}\left(p^1_{\epsilon_m}+p^2_{\epsilon_m}\right)
\end{split}
\qquad\mbox{in $\mathcal{D}'(\omega)$}.
\end{equation}

\section{The Lipschitz truncation}\label{Sect-Lip-T}

To start this section, let us set
\begin{equation}\label{v-em-chi}
\mathbf{w}_{\epsilon_m}:=(\mathbf{u}_{\epsilon_m}-\mathbf{u})\chi_{\omega},
\end{equation}
where $\chi_{\omega}$ denotes the characteristic function of the set $\omega$ introduced in (\ref{om-C2}). %
Having in mind the extension of (\ref{dist2-der}) to $\mathds{R}^{N}$, here we shall consider that
\begin{equation}\label{Ups}
\mathbf{\Upsilon}_{\epsilon_m}:=\mathbf{\Upsilon}^1_{\epsilon_m}+\mathbf{\Upsilon}^2_{\epsilon_m}
\end{equation}
is extended from $\omega$ to $\mathds{R}^{N}$ by zero, where
\begin{equation}\label{Ups-1}
\mathbf{\Upsilon}^1_{\epsilon_m}:=-\left(\mathbf{S}(\mathbf{D}(\mathbf{u}_{\epsilon_m}))-\mathbf{S}+\epsilon_m|\mathbf{D}(\mathbf{u}_{\epsilon_m})|^{\beta-2}\mathbf{D}(\mathbf{u}_{\epsilon_m})\right)+p^1_{\epsilon_m}\mathbf{I},
\end{equation}
\begin{equation}\label{Ups-2}
\mathbf{\Upsilon}^2_{\epsilon_m}:=\mathbf{u}_{\epsilon_m}\otimes\mathbf{u}_{\epsilon_m}\Phi_{\epsilon_m}(|\mathbf{u}_{\epsilon_m}|)-                                    \mathbf{u}\otimes\mathbf{u}+p^2_{\epsilon_m}\mathbf{I},
\end{equation}
and $\mathbf{I}$ denotes the identity tensor. %
Now, due to the definition (\ref{v-em-chi}) and by virtue of (\ref{convg-La}) and (\ref{convg-s-g}), we have
\begin{equation}\label{un-b-v-em}
\mbox{$\mathbf{w}_{\epsilon_m}\to \mathbf{0}$\quad weakly in $\mathbf{W}^{1,{\alpha}}(\mathds{R}^{N})$,\quad as $m\to\infty$,}
\end{equation}
\begin{equation}\label{str-v-em}
\mbox{$\mathbf{w}_{\epsilon_m}\to \mathbf{0}$\quad strongly in $\mathbf{L}^{\gamma}(\mathds{R}^{N})$,\quad as $m\to\infty$,
\quad
for any $\gamma:1\leq\gamma<\alpha^{\ast}$.
}
\end{equation}
Moreover, due to (\ref{convg-S-b'}), (\ref{conv-Re-0}) and (\ref{bd2-p1-em}) by one hand, and due to (\ref{convg-Phiuxu}), (\ref{limit-eq-G}) and (\ref{bd2-p2-em}) on the other, we have
\begin{equation}\label{un-b-T-p1-em}
\|\mathbf{\Upsilon}^1_{\epsilon_m}\|_{\mathbf{L}^{\beta'}(\mathds{R}^{N})}\leq C,
\end{equation}
\begin{equation}\label{un-bA-uxu-p2-em}
\|\mathbf{\Upsilon}^2_{\epsilon_m}\|_{\mathbf{L}^{\frac{\alpha^{\ast}}{2}}(\mathds{R}^{N})}\leq C.
\end{equation}
In addition to (\ref{un-bA-uxu-p2-em}), we see that, due to (\ref{convg-s-g}) and (\ref{bd2-p2-em}),
\begin{equation}\label{un-b-uxu-p2-em}
\mathbf{\Upsilon}^2_{\epsilon_m}\to 0\quad\mbox{strongly in $\mathbf{L}^{\frac{\gamma}{2}}(\mathds{R}^{N})$,\quad as $m\to\infty$,
\quad
for any $\gamma:1\leq\gamma<\alpha^{\ast}$.}
\end{equation}

\noindent Next, let us consider the Hardy-Littlewood maximal functions of $|\mathbf{w}_{\epsilon_m}|$ and $|\mathbf{\nabla}\mathbf{w}_{\epsilon_m}|$ defined by
\begin{eqnarray*}
  &\displaystyle& \mathcal{M}(|\mathbf{w}_{\epsilon_m}|)(\mathbf{x}):=\sup_{0<R<\infty}\frac{1}{\mathcal{L}_N(B_R(\mathbf{x}))}\int_{B_R(\mathbf{x})}|\mathbf{w}_{\epsilon_m}(\mathbf{y})|\,d\mathbf{y},\\
   &\displaystyle& \mathcal{M}(|\mathbf{\nabla}\mathbf{w}_{\epsilon_m}|)(\mathbf{x}):=\sup_{0<R<\infty}\frac{1}{\mathcal{L}_N(B_R(\mathbf{x}))}\int_{B_R(\mathbf{x})}|\mathbf{\nabla}\mathbf{w}_{\epsilon_m}(\mathbf{y})|\,d\mathbf{y};
\end{eqnarray*}
where $B_R(\mathbf{x})$ denotes the ball of $\mathds{R}^N$ centered at $\mathbf{x}$ and with radius $R>0$, and $\mathcal{L}_N(\omega)$ is the $N$-dimensional Lebesgue measure of $\omega$. %
Arguing as in \cite[p. 218]{DMS-2008} and using the boundedness of the Hardy-Littlewood maximal operator $\mathcal{M}$, we can prove that for all $m\in\mathds{N}$ and all $j\in\mathds{N}_0$ there exists
\begin{equation}\label{int-2j}
\lambda_{m,j}\in\left[2^{2^j},2^{2^{j+1}}\right)
\end{equation}
such that
\begin{equation}\label{Leb-F}
\mathcal{L}_{N}\left(F_{m,j}\right)\leq
2^{-j}\lambda_{m,j}^{-\gamma}\,\|\mathbf{w}_{\epsilon_m}\|_{\mathbf{L}^{\gamma}(\mathds{R}^{N})},
\quad
\mbox{for any $\gamma:1\leq\gamma<\alpha^{\ast}$,}
\end{equation}
\begin{equation}\label{Leb-G}
\mathcal{L}_{N}\left(G_{m,j}\right)\leq
2^{-j}\lambda_{m,j}^{-\alpha}\,\|\mathbf{\nabla}\mathbf{w}_{\epsilon_m}\|_{\mathbf{L}^{\alpha}(\mathds{R}^{N})},
\end{equation}
where
\begin{eqnarray*}
  &\displaystyle& F_{m,j}:=\left\{\mathbf{x}\in\mathds{R}^{N}:\mathcal{M}(|\mathbf{w}_{\epsilon_m}|)(\mathbf{x})>2\lambda_{m,j}\right\},\\
   &\displaystyle& G_{m,j}:=\left\{\mathbf{x}\in\mathds{R}^{N}:\mathcal{M}(|\mathbf{\nabla}\mathbf{w}_{\epsilon_m}|)(\mathbf{x})>2\lambda_{m,j}\right\}.
\end{eqnarray*}
Setting
\begin{equation}\label{R-mj}
R_{m,j}:=F_{m,j}\cup G_{m,j}\cup\left\{\mathbf{x}\in\mathds{R}^{N}:\ \mbox{$\mathbf{x}$ is not a Lebesgue point of $|\mathbf{w}_{\epsilon_m}|$}\right\},
\end{equation}
we can see that, by virtue of (\ref{Leb-F})-(\ref{R-mj}),
\begin{equation}\label{Leb-R}
\mathcal{L}_{N}\left(R_{m,j}\right)\leq
2^{-j}\lambda_{m,j}^{-\alpha}\,\|\mathbf{w}_{\epsilon_m}\|_{\mathbf{W}^{1,{\alpha}}(\mathds{R}^{N})}.
\end{equation}
In addition, due to (\ref{un-b-v-em})-(\ref{str-v-em}) and (\ref{int-2j}),
\begin{equation}\label{lim-m-R}
\limsup_{m\to\infty}\mathcal{L}_{N}\left(R_{m,j}\right)\leq C2^{-j}\lambda_{m,j}^{-\alpha}.
\end{equation}

\noindent Then, by \cite{AF-1988} together with (\ref{v-em-chi}),
\begin{equation}\label{z-mj}
\exists\ \mathbf{z}_{m,j}\in\mathbf{W}^{1,\infty}(\mathds{R}^N),\qquad
\mathbf{z}_{m,j}=\left\{\begin{array}{cc}
                   \mathbf{w}_{\epsilon_m} & \mbox{in $\omega\setminus A_{m,j}$} \\
                   0 & \mathds{R}^N\setminus\omega
                 \end{array}\right.\,,
\end{equation}
where
\begin{equation}\label{A-mj}
A_{m,j}:=\{\mathbf{x}\in\omega:\mathbf{z}_{m,j}(\mathbf{x})\not=\mathbf{w}_{\epsilon_m}(\mathbf{x})\},
\end{equation}
such that
\begin{equation}\label{eq1-AF}
\|\mathbf{z}_{m,j}\|_{\mathbf{L}^{\infty}(\omega)}\leq2\lambda_{m,j},
\end{equation}
\begin{equation}\label{eq2-AF}
\|\mathbf{\nabla}\mathbf{z}_{m,j}\|_{\mathbf{L}^{\infty}(\omega)}\leq C\lambda_{m,j},\quad C=C(N,\omega).
\end{equation}
Moreover, by  \cite[Proposition 2.2]{Landes-1996} and using (\ref{Leb-F})-(\ref{R-mj}) and (\ref{A-mj}),
\begin{equation}\label{eq3-AF}
A_{m,j}\subset\omega\cap R_{m,j}.
\end{equation}
As a consequence of (\ref{eq3-AF}) together with (\ref{Leb-R}) and (\ref{lim-m-R}),
\begin{equation}\label{Leb-A}
\mathcal{L}_{N}\left(A_{m,j}\right)\leq
2^{-j}\lambda_{m,j}^{-\alpha}\,\|\mathbf{w}_{\epsilon_m}\|_{\mathbf{W}^{1,{\alpha}}(\mathds{R}^{N})},
\end{equation}
\begin{equation}\label{lim-m-A}
\limsup_{m\to\infty}\mathcal{L}_{N}\left(A_{m,j}\right)\leq C
2^{-j}\lambda_{m,j}^{-\alpha}.
\end{equation}
On the other hand, due to (\ref{eq1-AF})-(\ref{eq2-AF}) and (\ref{un-b-v-em}) together with (\ref{Leb-A})-(\ref{lim-m-A}), we can prove that for any $j\in\mathds{N}_0$
\begin{equation}\label{w-z-em-a}
\mbox{$\mathbf{z}_{m,j}\to \mathbf{0}$\quad weakly in $\mathbf{W}^{1,{\alpha}}_0(\omega)$,\quad as $m\to\infty$.}
\end{equation}
Then by Sobolev's compact imbedding theorem, we get for any $j\in\mathds{N}_0$
\begin{equation*}
\mbox{$\mathbf{z}_{m,j}\to \mathbf{0}$ strongly in $\mathbf{L}^{\gamma}(\omega)$,\quad as $m\to\infty$,\quad for any $\gamma: 1\leq\gamma<\alpha^{\ast}$.}
\end{equation*}
Using this information, (\ref{eq1-AF}) and interpolation, we prove that for any $j\in\mathds{N}_0$
\begin{equation}\label{str-z-em}
\mbox{$\mathbf{z}_{m,j}\to \mathbf{0}$\quad strongly in $\mathbf{L}^s(\omega)$,\quad as $m\to\infty$,\quad for any $s:1\leq s<\infty$.}
\end{equation}
Finally, as a consequence of (\ref{w-z-em-a}) and (\ref{str-z-em}), we obtain for any $j\in\mathds{N}_0$
\begin{equation}\label{w-con-s-z-em}
\mbox{$\mathbf{z}_{m,j}\to \mathbf{0}$\quad weakly in $\mathbf{W}^{1,s}_0(\omega)$,\quad as $m\to\infty$,\quad for any $s:1\leq s<\infty$.}
\end{equation}

\section{Convergence of the approximated extra stress tensor}\label{Sect-Conv-EST}

Let us first observe that, using the notations (\ref{Ups})-(\ref{Ups-2}), we can write (\ref{dist2-der}) as
\begin{equation}\label{dist3-der}
\mathbf{div}\mathbf{\Upsilon}_{\epsilon_m}=\mathbf{0}
\quad\mbox{in\quad $\mathcal{D}'(\omega)$}.
\end{equation}
On the other hand, due to (\ref{un-b-T-p1-em})-(\ref{un-bA-uxu-p2-em}), $\mathbf{\Upsilon}_{\epsilon_m}\in\mathbf{L}^{r}(\mathds{R}^N)$ for $r$ satisfying to (\ref{Q-in-Lr}). %
Then, using this information and (\ref{w-con-s-z-em}), we infer, from (\ref{dist3-der}), that for any $j\in\mathds{N}_0$
\begin{equation}\label{dist4-der}
  \int_{\omega}\mathbf{\Upsilon}_{\epsilon_m}:\mathbf{\nabla}\mathbf{z}_{m,j}\,d\mathbf{x}=0.
\end{equation}
Expanding $\mathbf{\Upsilon}_{\epsilon_m}$ in (\ref{dist4-der}) through the notations (\ref{Ups})-(\ref{Ups-2}) and subtracting and adding the integral of $\mathbf{S}(\mathbf{D}(\mathbf{u})):\mathbf{D}(\mathbf{z}_{m,j})$ to the left hand side of the resulting equation, we obtain for any $j\in\mathds{N}_0$
\begin{equation}\label{dist5-der}
\begin{split}
  & \int_{\omega}\left(\mathbf{S}(\mathbf{D}(\mathbf{u}_{\epsilon_m}))-\mathbf{S}(\mathbf{D}(\mathbf{u}))\right):\mathbf{D}(\mathbf{z}_{m,j})\,d\mathbf{x} =\int_{\omega}\left(\mathbf{S}-\mathbf{S}(\mathbf{D}(\mathbf{u}))\right):\mathbf{D}(\mathbf{z}_{m,j})\,d\mathbf{x}\\
  &-\int_{\omega}\epsilon_m|\mathbf{D}(\mathbf{u}_{\epsilon_m})|^{\beta-2}\mathbf{D}(\mathbf{u}_{\epsilon_m}):\mathbf{D}(\mathbf{z}_{m,j})\,d\mathbf{x}+
  \int_{\omega}p^1_{\epsilon_m}\,\mathrm{div}\,\mathbf{z}_{m,j}\,d\mathbf{x}+\\
  &\int_{\omega}\left(\mathbf{u}_{\epsilon_m}\otimes\mathbf{u}_{\epsilon_m}\Phi_{\epsilon_m}(|\mathbf{u}_{\epsilon_m}|)-                                    \mathbf{u}\otimes\mathbf{u}+p^2_{\epsilon_m}\mathbf{I}\right):\mathbf{D}(\mathbf{z}_{m,j})\,d\mathbf{x}\\
  &:=J_{m,j}^1+J_{m,j}^2+J_{m,j}^3+J_{m,j}^4.
  \end{split}
\end{equation}
We claim that, for a fixed $j$,
\begin{equation}\label{dist6-der}
\lim_{m\to\infty}\int_{\omega}\left(\mathbf{S}(\mathbf{D}(\mathbf{u}_{\epsilon_m}))-\mathbf{S}(\mathbf{D}(\mathbf{u}))\right):\mathbf{D}(\mathbf{z}_{m,j})\,d\mathbf{x}
\leq C2^{-\frac{j}{\beta}}
\end{equation}
To prove this, we will carry out the passage to the limit $m\to\infty$ in all absolute values $|J_{m,j}^i|$, $i=1,\dots,4$.

$\bullet\ \limsup_{m\to\infty}|J_{m,j}^1|=0$. By (\ref{w-con-s-z-em}), with $s=\beta$, this is true once we can justify that
$\mathbf{S}-\mathbf{S}(\mathbf{D}(\mathbf{u}))$ is uniformly bounded in $\mathbf{L}^{\beta'}(\omega)$. %
But this is a consequence of (\ref{convg-S-b'}), the continuous imbedding $\mathbf{L}^{q'(\cdot)}(\omega)\hookrightarrow\mathbf{L}^{\beta'}(\omega)$ and (\ref{est-q'}).

$\bullet\ \limsup_{m\to\infty}|J_{m,j}^2|=0$. Indeed, by Hölder's inequality, (\ref{eq2-AF}), (\ref{int-2j}) and (\ref{e-inequality-qr}),
we have successively
\begin{equation*}
\begin{split}
  |J_{m,j}^2|&\leq \|\epsilon_m|\mathbf{D}(\mathbf{u}_{\epsilon_m})|^{\beta-2}\mathbf{D}(\mathbf{u}_{\epsilon_m})\|_{\mathbf{L}^{1}(\omega)}
\|\mathbf{\nabla}\mathbf{z}_{m,j}\|_{\mathbf{L}^{\infty}(\omega)} \\
       &\leq C_1\lambda_{m,j}\, \epsilon_m^{\frac{1}{\beta}}\left(\int_{\omega}\epsilon_m|\mathbf{D}(\mathbf{u}_{\epsilon_m})|^{\beta}\,d\mathbf{x}\right)^{\frac{\beta-1}{\beta}}\leq C_2\,\epsilon_m^{\frac{1}{\beta}}\to 0,\quad\mbox{as}\ m\to\infty.
\end{split}
\end{equation*}

$\bullet\ \limsup_{m\to\infty}|J_{m,j}^3|\leq C 2^{-\frac{j}{\beta}}$. In fact, by Hölder's inequality and (\ref{bd2-p1-em}) together with (\ref{convg-S-b'}) and (\ref{conv-Re-0}), and using (\ref{z-mj})) together with (\ref{v-em-chi}),
\begin{equation*}
|J_{m,j}^3|\leq C_1\|\mathrm{div}\,\mathbf{z}_{m,j}\|_{\mathbf{L}^{\beta}(\omega)}
     \leq C_1\|\mathbf{\nabla}\mathbf{z}_{m,j}\|_{\mathbf{L}^{\infty}(\omega)}\mathcal{L}_N(A_{m,j})^{\frac{1}{\beta}}.
\end{equation*}
The result follows by the application of (\ref{eq2-AF}), (\ref{int-2j}) and (\ref{lim-m-A}), provided that
\begin{equation}\label{cond1-HI}
\lambda_{m,j}^{1-\frac{\alpha}{\beta}}\quad\mbox{is uniformly bounded in $m$.}
\end{equation}

$\bullet\ \limsup_{m\to\infty}|J_{m,j}^4|=0$. Using Hölder's inequality and (\ref{Ups-2}), we have
\begin{equation*}
  |J_{m,j}^4|\leq \|\mathbf{\Upsilon}^2_{\epsilon_m}\|_{\mathbf{L}^{1}(\omega)}\|\mathbf{\nabla}\mathbf{z}_{m,j}\|_{\mathbf{L}^{\infty}(\omega)}\leq C_1\|\mathbf{\Upsilon}^2_{\epsilon_m}\|_{\mathbf{L}^{1}(\omega)}\to 0,\quad\mbox{as}\ m\to\infty.
\end{equation*}
The last inequality and the conclusion follow, respectively, from (\ref{eq2-AF}) and (\ref{int-2j}), and (\ref{un-b-uxu-p2-em}) with $\gamma=2$, observing here that assumption (\ref{al-TE}) implies $2<\alpha^{\ast}$.

Gathering the estimates above we just have proven (\ref{dist6-der}).

We proceed with the proof by using an argument due to \cite[Theorem 5]{DalM-M}.
Firstly, observing the definition of $\mathbf{z}_{m,j}$ (\emph{cf.} (\ref{z-mj})), we have
\begin{equation}\label{eq21-est-theta}
\int_{\omega}\left(\mathbf{S}(\mathbf{D}(\mathbf{u}_{\epsilon_m}))-\mathbf{S}(\mathbf{D}(\mathbf{u}))\right):\mathbf{D}(\mathbf{z}_{m,j})\,d\mathbf{x}= I_{m,j}+II_{m,j},\\ \end{equation}
where
\begin{equation*}
\begin{split}
&I_{m,j}:=\int_{\omega\setminus A_{m,j}}\left(\mathbf{S}(\mathbf{D}(\mathbf{u}_{\epsilon_m}))-\mathbf{S}(\mathbf{D}(\mathbf{u}))\right):(\mathbf{D}(\mathbf{u}_{\epsilon_m})-\mathbf{D}(\mathbf{u}))\,d\mathbf{x},\\
&II_{m,j}:= \int_{A_{m,j}}\left(\mathbf{S}(\mathbf{D}(\mathbf{u}_{\epsilon_m}))-\mathbf{S}(\mathbf{D}(\mathbf{u}))\right):\mathbf{D}(\mathbf{z}_{m,j})\,d\mathbf{x}.
     \end{split}
\end{equation*}
Then (\ref{dist6-der}) and (\ref{eq21-est-theta}) imply that
\begin{equation}\label{est-I1}
\limsup_{m\to\infty}I_{m,j}\leq\limsup_{m\to\infty}|II_{m,j}|+C2^{-\frac{j}{\beta}}.
\end{equation}
For the term $II_{m,j}$, we have by applying successively Hölder's inequality, (\ref{convg-S-b'}), the continuous imbedding $\mathbf{L}^{q'(\cdot)}(\omega)\hookrightarrow\mathbf{L}^{\beta'}(\omega)$ and (\ref{est-q'}) altogether with (\ref{eq2-AF}),
\begin{equation*}
\begin{split}
|II_{m,j}|&\leq
C_1
\|\mathbf{S}(\mathbf{D}(\mathbf{u}_{\epsilon_m}))-\mathbf{S}(\mathbf{D}(\mathbf{u}))\|_{\mathbf{L}^{\beta'}(A_{m,j})}
\|\mathbf{\nabla}\mathbf{z}_{m,j}\|_{\mathbf{L}^{\beta}(A_{m,j})}\\
&\leq
C_2\lambda_{m,j}\mathcal{L}_N(A_{m,j})^{\frac{1}{\beta}}.
\end{split}
\end{equation*}
Then, (\ref{int-2j}) and (\ref{lim-m-A}) yield that for  any $j\in\mathds{N}_0$
\begin{equation}\label{est-I2}
\limsup_{m\to\infty}|II_{m,j}|\leq
C2^{-\frac{j}{\beta}}\lambda_{m,j}^{1-\frac{\alpha}{\beta}}.
\end{equation}
As a consequence of (\ref{est-I1}) and (\ref{est-I2}), we obtain for any $j\in\mathds{N}_0$
\begin{equation}\label{est2-I1}
\limsup_{m\to\infty}|I_{m,j}|\leq
C2^{-\frac{j}{\beta}}\left(1+\lambda_{m,j}^{1-\frac{\alpha}{\beta}}\right).
\end{equation}
Arguing as we did to prove (\ref{est-I2})-(\ref{est2-I1}) and using (\ref{z-mj}) and (\ref{lim-m-A}), we have for any $\theta\in(0,1)$
\begin{equation}\label{eq3-est-theta}
\limsup_{m\to\infty}\int_{\omega}g_{\epsilon_m}^{\theta}\,d\mathbf{x}
\leq
C_1 2^{-\theta\frac{j}{\beta}}\left(1+\lambda_{m,j}^{1-\frac{\alpha}{\beta}}\right)^{\theta}+C_2 2^{-\theta\frac{j}{\beta}-(1-\theta)j}\lambda_{m,j}^{\left(1-\frac{\alpha}{\beta}\right)\theta-(1-\theta)\alpha}\,,
\end{equation}
where
\begin{equation*}
g_{\epsilon_m}:=\left(\mathbf{S}(\mathbf{D}(\mathbf{u}_{\epsilon_m}))-\mathbf{S}(\mathbf{D}(\mathbf{u}))\right):(\mathbf{D}(\mathbf{u}_{\epsilon_m})-\mathbf{D}(\mathbf{u})).
\end{equation*}
Since $\beta>1$, $\theta\in(0,1)$ and $j\in\mathds{N}_0$ is arbitrary, $2^{-\theta\frac{j}{\beta}}\to 0$ and $2^{-\theta\frac{j}{\beta}-(1-\theta)j}\to 0$, as $j\to\infty$.
This and (\ref{eq3-est-theta}) imply that for any $\theta\in(0,1)$
\begin{equation*}
\limsup_{m\to\infty}\int_{\omega}g_{\epsilon_m}^{\theta}\,d\mathbf{x}=0
\end{equation*}
provided that (\ref{cond1-HI}) holds. %
Then, passing to a subsequence,
\begin{equation}\label{eq2-est-theta}
g_{\epsilon_m}\to 0\quad\mbox{a.e. in $\omega$},\quad\mbox{as $m\to\infty$}.
\end{equation}
From the continuity of $\mathbf{S}$ on $\mathbf{D}(\mathbf{u})$ (\emph{cf.} condition (A)), the strict monotonicity condition (D), (\ref{eq2-est-theta}) and \cite[Lemma 6]{DalM-M} (see also \cite[Lemme 2.2.2]{Lions-1969}),
\begin{equation}\label{eq4-est-theta}
\mathbf{D}(\mathbf{u}_{\epsilon_m})\to\mathbf{D}(\mathbf{u})\quad\mbox{a.e. in $\omega$},\quad\mbox{as $m\to\infty$}.
\end{equation}
Finally, (\ref{est-gam'}) and (\ref{eq4-est-theta}) allow us to use Vitali's theorem together with (\ref{convg-S-b'}) to conclude that $\mathbf{S}=\mathbf{S}(\mathbf{D}(\mathbf{u}))$.

\section{Answer to Question~\ref{th-exst-ws-pp}}

From Section~\ref{Sect-Exist-RP} until Section~\ref{Sect-Conv-EST} we have proven the existence of, at least, a weak solution to the problem (\ref{geq1-inc})-(\ref{geq1-bc}) in the sense of Definition~\ref{weak-sol-vq} and satisfying to the conditions stated in Question~\ref{th-exst-ws-pp}, provided condition (\ref{cond1-HI}) is fulfilled. %
A simple analysis shows us that condition (\ref{cond1-HI}) is equivalent to assume that $\alpha\geq\beta$. But this cannot happen unless  $\alpha=\beta$. In this situation, we would fall in the case of a constant exponent $q$ studied in \cite{FMS-2003}.
If we go further behind, we see that condition (\ref{cond1-HI}) came as a result of (\ref{lim-m-A}) and this in turn had its origin in (\ref{un-b-v-em})-(\ref{str-v-em}). Therefore the best way to assure that (\ref{cond1-HI}) is fulfilled is to assume that
(\ref{un-b-v-em})-(\ref{str-v-em}) are satisfied with $\alpha$ replaced by $\beta$, \emph{i.e.}
\begin{equation}\label{un-b-v-em-beta}
\mbox{$\mathbf{w}_{\epsilon_m}\to \mathbf{0}$\quad weakly in $\mathbf{W}^{1,{\beta}}(\mathds{R}^{N})$,\quad as $m\to\infty$,}
\end{equation}
\begin{equation}\label{str-v-em-beta}
\mbox{$\mathbf{w}_{\epsilon_m}\to \mathbf{0}$\quad strongly in $\mathbf{L}^{\gamma}(\mathds{R}^{N})$,\quad as $m\to\infty$,
\quad
for any $\gamma:1\leq\gamma<\beta^{\ast}$.
}
\end{equation}
We observe that (\ref{un-b-v-em})-(\ref{str-v-em}) came as a result of (\ref{est-inf-gam}).
In consequence, (\ref{un-b-v-em-beta})-(\ref{str-v-em-beta}) hold if we had
\begin{equation}\label{est-inf-beta}
\|\mathbf{u}_{\epsilon}\|_{\mathbf{V}_{\beta}}\leq C
\end{equation}
instead of (\ref{est-inf-gam}). %
Finally, condition (\ref{est-inf-beta}) is satisfied if
\begin{equation}\label{cond2-HI}
\int_{\Omega}|\mathbf{D}(\mathbf{u}_{\epsilon})|^{\beta}\,d\mathbf{x}\leq C\,.
\end{equation}
In consequence, if we go back to Sections~\ref{Sect-Lip-T} and~\ref{Sect-Conv-EST} and replace all the exponents $\alpha$ by $\beta$, than
 we have $\lambda_{m,j}^{1-\frac{\alpha}{\beta}}=1$ in condition (\ref{cond1-HI}), and the existence result follows. %
Condition (\ref{cond2-HI}) can be seen as a consequence of the following higher integrability condition: assume that exists $\delta>0$ such that
\begin{equation}\label{HI-cond}
\int_{\omega'}|\mathbf{D}(\mathbf{u})|^{q(\mathbf{x})(1+\delta)}d\mathbf{x}<\infty
\end{equation}
for any subdomain $\omega'\subset\subset\Omega$.
Then, under all the assumptions of Question~\ref{th-exst-ws-pp}, we can prove the existence of weak solutions for the problem (\ref{geq1-inc})-(\ref{geq1-bc}) with $q=q(\mathbf{x})$, which \emph{a priori} satisfy condition (\ref{HI-cond}).
This property is crucial to control the gradients of velocity in the space $\mathbf{L}^{\beta}(\Omega)$ and with the technique we used we can control them only in the spaces $\mathbf{L}^{\alpha}(\Omega)$ and $\mathbf{L}^{q(\cdot)}(\Omega)$. %
Despite assumption (\ref{HI-cond}) is so strong that weakens very much such an existence result, we observe that this property is satisfied by the weak solutions to some fluid problems. In fact, fluids with viscosity dependence described using non-standard growth conditions have been treated, in the stationary case, in various settings (see \cite{AMS-2004,AZ-2005} and the references therein).
We think that the approach followed in \cite{AMS-2004, AZ-2005} can be potentially
useful to extend theses results to our problem. %
Therefore we are let to believe the higher integrability property expressed by assumption (\ref{HI-cond}) is satisfied by
every weak solution to the problem (\ref{geq1-inc})-(\ref{geq1-bc}) with $q=q(\mathbf{x})$. %
In this case, to prove an existence result for our problem we do not need the log-Hölder continuity condition (\ref{glob-log-H}) on $q$. %
On the other hand, we can realize that for models of generalized fluid flows in which the stress tensor brings itself this higher regularity, the existence result follows without assuming (\ref{glob-log-H}) and (\ref{HI-cond}). An example of this situation is the problem with the stress tensor defined by
$$\mathbf{S}=\left(\mu+\tau|\mathbf{D}|^{q(\mathbf{x})-2}\right)\mathbf{D}\,,$$
where $\mu$ and $\tau$ are positive constants related with the viscosity. In this case, an existence result can be proved for $\frac{2N}{N+2}<\alpha<\beta\leq 2$ proceeding as in the above sections.

A similar analysis can be done for the parabolic version of the problem (\ref{geq1-inc})-(\ref{geq1-bc}). %
In this case all the reasoning is identical and we just have to use the parabolic versions of the results considered from Section~\ref{Sect-Exist-RP} to Section~\ref{Sect-Conv-EST}. %
The impact of our work in the transient problem is in fact more important, because the parabolic extension of the work \cite{DMS-2008} is, to the best of our knowledge, still not proved. %
A thorough analysis of these problems is being written and it will be published elsewhere.

\appendix

\end{document}